\documentclass[oneside, a4paper,11pt,reqno]{amsart}
\providecommand{\U}[1]{\protect\rule{.1in}{.1in}}
\newcommand\bkE{{\mathbb {E}}}
\newcommand\E{{\mathbb {E}}}

\newtheorem {Lemma}{Lemma}[section]
\newtheorem {Claim}{Claim}[section]
\newtheorem {Theorem}{Theorem}[section]
\newtheorem {Proposition}{Proposition}[section]

\newcommand\beq{\begin{equation}}
\newcommand\eeq{\end{equation}}

\def\esup{\mathop{\bf esssup}}

\begin{document}
\begin{center} {\bf \Large Rates of convergence  in
the strong invariance principle for non adapted sequences.
Application to ergodic automorphisms of the torus}\vskip15pt

J\'er\^ome Dedecker $^{a}$, Florence Merlev\`{e}de $^{b}$ {\it
and\/} Fran\c{c}oise P\`ene  $^{c}$
\end{center}
\vskip10pt
$^a$ Universit\'e Paris Descartes, Sorbonne Paris Cit\'e, Laboratoire MAP5
and CNRS UMR 8145. Email: jerome.dedecker@parisdescartes.fr\\ \\
$^b$ Universit\'e Paris Est, LAMA and CNRS UMR 8050.\\
E-mail: florence.merlevede@univ-mlv.fr\\ \\
$^c$ Universit\'e de Brest,  Laboratoire de Math\'ematiques de Bretagne
Atlantique UMR CNRS 6205.
E-mail: francoise.pene@univ-brest.fr\vskip10pt
{\it Key words}: almost sure invariance principle, strong approximations, non adapted sequences, ergodic automorphisms of the torus.\vskip5pt


{F. P\`ene is partially supported by the french ANR projects MEMEMO2 and PERTURBATIONS.}
\begin{center}
{\bf Abstract}\vskip10pt
\end{center}

In this paper, we give rates of convergence in the strong invariance principle
for non-adapted sequences satisfying projective
criteria. The results apply to the iterates  of ergodic automorphisms
$T$ of the $d$-dimensional torus ${\mathbb T}^d$, even in the non hyperbolic case.
In this context,
we give a large class of unbounded function
$f$ from ${\mathbb T}^d$ to ${\mathbb R}$, for which the partial sum
$f\circ T+ f \circ T^2 + \cdots + f \circ T^n$ satisfies a strong invariance
principle with an explicit rate of convergence.

\section{Introduction and notations}
\setcounter{equation}{0}

\setcounter{equation}{0}

Let $(\Omega,{\mathcal A}, {\mathbb P} )$ be a
probability space, and $T:\Omega \mapsto \Omega$ be
 a bijective bimeasurable transformation preserving the probability ${\mathbb P} $.
For a $\sigma$-algebra ${\mathcal F}_0 $ satisfying ${\mathcal F}_0
\subseteq T^{-1 }({\mathcal F}_0)$, we define the nondecreasing
filtration $({\mathcal F}_i)_{i \in {\mathbb Z}}$ by ${\mathcal F}_i =T^{-i
}({\mathcal F}_0)$.
The ${\mathbb L}^{p} $ norm of a random variable $X$ is denoted by
$\|X \|_{p}= (  \E ( | X|^p)  )^{1/p}$.

Let $X_0$ be a real-valued and square integrable random variable such that
${\mathbb E}(X_0)=0$, and
define the stationary sequence  $(X_i)_{i \in \mathbb Z}$ by
$X_i = X_0 \circ T^i$. Define then the partial sum by $S_n= X_1+X_2+ \cdots  + X_n$.
According
to the Birkhoff-Khinchine theorem, $S_n$ satisfies a strong law of large numbers.
One can go further in the study of the statistical properties of $S_n$.
We study here the rate of convergence in the almost sure
invariance principle (ASIP).
More precisely,
 we give  conditions under which there exists a sequence
of independent identically distributed (iid) Gaussian random variables $(Z_i)_{i\geq 1}$ such that
\begin{equation}\label{ASIPp}
  \sup_{1 \leq k \leq n} \Big | \sum_{i=1}^k (X_i-Z_i) \Big|= o(n^{1/p}L(n)) \quad
  \text{almost surely,}
\end{equation}
for $p\in ]2,4]$ and $L$ an explicit slowly varying function. Let us recall that, in the iid case, Koml\'os,
Major and Tusn\'ady \cite{KMT} and Major \cite{Major} obtained an ASIP with the optimal rate $o(n^{1/p})$ in \eqref{ASIPp} as soon as the random variables admit a moment of
order $p$.

Since the seminal paper by  Philipp and Stout \cite{PS}, many authors have considered
this problem in a dependent context, but most of the papers
deal with the adapted case, when  $X_0$ is ${\mathcal F}_0$ measurable (for instance,
${\mathcal F}_0$ is the past $\sigma$-algebra $\sigma( X_i, i \leq 0)$).
Unfortunately, it is quite common to encounter dynamical systems for which the natural filtration does not allow to control any quantity involving terms of the type
$\|{\mathbb E}(X_n|{\mathcal F}_0)\|_p$.

In this paper, we shall not assume that $X_0$ is ${\mathcal F}_0$-measurable,
and we shall give conditions on
$\|{\mathbb E}(X_n|{\mathcal F}_0)\|_p$, $\|X_{-n}-{\mathbb E}(X_{-n}|{\mathcal F}_0)\|_p$
and
$\|{\mathbb E}(S_n^2|{\mathcal F}_{-n})-{\mathbb E}(S_n^2)\|_{p/2}$ for
(\ref{ASIPp}) to hold (see Theorems \ref{ThNAS} and \ref{ThNASp=4} of Section \ref{MR}).
These conditions are in the same spirit as
those given by  Gordin \cite{Go} for $p=2$ to get  the usual central limit theorem.
Our proof is based on the approximation
$$\sum_{i=1}^nX_i=M_n+R_n$$
by the martingale
$M_n= d_1+d_2+ \cdots +d_n$, where $d_i$ is the martingale
difference
$$
d_i=\sum_{k \in {\mathbb Z}} ({\mathbb E}(X_k|{\mathcal F}_i)-
{\mathbb E}(X_k|{\mathcal F}_{i-1}))
$$
introduced by Gordin \cite{Go} and Heyde \cite{Heyde}.
In the adapted case, similar conditions are
given in the recent paper \cite{DDM},
together with a long list of applications.

In the non adapted case, it is easy to see that our results apply
to  a large class of two-sided functions
of iid sequences, or two-sided functions of absolutely regular sequences.
But they also apply to much complicated dynamical systems, for which
such a representation by functions of absolutely regular sequences is
not available. In the next section, we consider the case where $T$ is
an ergodic automorphism of
 the $d$-dimensional torus $\mathbb T^d$, and ${\mathbb P}$ is the Lebesgue measure on
 $\mathbb T^d$. In this context, we use the
 $\sigma$-algebra ${\mathcal F}_i$ considered by Le Borgne \cite{SLB}.
As a consequence of Theorem \ref{auto-tore2}, we obtain that
(\ref{ASIPp}) holds for $p=4$  and $X_i=f\circ T^i$, where $f:{\mathbb T}^d \rightarrow {\mathbb R}$, as soon as
the Fourier coefficients $(c_{\mathbf k})_{\mathbf k\in\mathbb Z^d}$ of $f$ are such that
$$
  |c_{\mathbf k}| \leq A \prod_{i=1}^d \frac{1}{(1+|k_i|)^{3/4}
  \log^{\alpha}(2+|k_i|)} \quad \text{for some $\alpha>13/8$.}
$$
We also get that there exists a positive $\varepsilon$ such that
$$
\sup_{1 \leq k \leq n} \Big | \sum_{i=1}^k (X_i-Z_i) \Big|= o(n^{1/2-\varepsilon}) \quad
  \text{almost surely,}
$$
as soon as
$$
|c_{\mathbf k}| \leq A \prod_{i=1}^d \frac{1}{(1+|k_i|)^{\delta}}
\quad \text{for some $\delta>1/2$.}
$$
These rates of convergence in the almost sure invariance principle
complement the results by Leonov \cite{Leonov} and Le Borgne \cite{SLB}
for the central limit theorem and the almost sure invariance principle respectively.
Let us mention that Dolgopyat \cite{Dolgopyat} established an ASIP with the rate
$o(n^{1/2-\varepsilon})$ (for some $\varepsilon>0$) valid for
ergodic automorphisms of the torus and $f$ a H\"older continuous function.
Thanks to the decorrelation estimates obtained in
\cite{SLBFP}, 
the rate for H\"older observables can be improved
by applying the general result of Gou\"ezel in \cite{Gouezel}
to get the rate  $o(n^{1/4+\varepsilon})$ for every $\varepsilon>0$,
and by applying the results of the present paper to get the rate $o(n^{1/4}L(n))$.
Up to our knowledge, the present work gives the first strong approximations results for such
partially hyperbolic transformations $T$ for unbounded (and then non
continuous) functions  $f$.

To conclude, let us mention some previous works in the context of dynamical systems:
several results have been established with the rate  $o(n^{1/2-\varepsilon})$
for some $\varepsilon>0$
(see \cite{HofbauerKeller,DenkerPhilipp,Dolgopyat,Nagayama,MelbourneNicol1}).
Results giving a rate in $o(n^{1/4+\varepsilon})$
for every $\varepsilon>0$ can be found in
\cite{MelbourneTorok,FieldMelbourneTorok,MelbourneNicol2,Gouezel}. Most of
these results hold for bounded functions $f$.

Let us precise once again that we can reach the rate $o(n^{1/4}L(n))$
instead of $o(n^{ 1/4+\varepsilon})$
for every $\varepsilon>0$.
Moreover, our conditions giving the rate $o(n^{1/p}L(n))$ are related to
moments of order $p$ of $f$. Such results are not very common in the context of dynamical systems
(let us mention \cite{Gouezel} in the particular case of Gibbs-Markov maps, and 
\cite{DGM,MR} for generalized Pommeau-Manneville maps).

\section{ASIP with rates for ergodic automorphisms of the torus}

Let $d\ge 2$. We consider a group automorphism $T$ of the torus $\mathbb T^d=\mathbb R^d/\mathbb Z^d$.
For every $x\in\mathbb R^d$, we write $\bar x$ its class in $\mathbb T^d$.
We recall that $T$ is the quotient map of a linear map
$\tilde T:\mathbb R^d\rightarrow\mathbb R^d$
given by $\tilde T(x)=S\cdot x$, where $S$ is a $d\times d$-matrix with integer entries and with
determinant 1 or -1. The map $x\mapsto S\cdot x$
preserves the infinite Lebesgue measure $\lambda$ on $\mathbb R^d$
and $T$ preserves the probability Lebesgue measure $\bar\lambda$.
We suppose $T$ ergodic, which is equivalent to the fact that no eigenvalue of $S$
is a root of the unity. In this case, it is known that the spectral radius of $S$
is larger than one (and so $S$ admits at least an eigenvalue of modulus
larger than one and at least an eigenvalue of modulus smaller than one).
This hypothesis holds true in the case of hyperbolic automorphisms of
the torus (i.e. in the case when no  eigenvalue of $S$ has modulus one) but is much weaker.
Indeed, as mentioned in \cite{SLB}, the following matrix gives an example of an ergodic
non hyperbolic automorphism of $\mathbb T^4$~:
$$S:=\left(\begin{array}{cccc}0&0&0&-1\\1&0&0&2\\0&1&0&0\\0&0&1&2\end{array}\right).$$
When $T$ is ergodic and  non hyperbolic, the dynamical system $(\mathbb T^d,T,\bar\lambda)$ has no
Markov partition. However, it is possible to construct some measurable
partition \cite{Lind}, to prove a central limit theorem \cite{Leonov}. Moreover, in \cite{SLB}, Le Borgne
proved the functional central limit theorem and the Strassen strong invariance principle
for $(X_k=f\circ T^k)_k$ under weak hypotheses on $f$, thanks to Gordin's method and to the
partitions studied by Lind in \cite{Lind}.

We give here rates of convergence in the strong invariance
principle  for $(X_k=f\circ T^k)_k$
under conditions on the Fourier coefficients of $f:{\mathbb T}^d \rightarrow {\mathbb R}$. In what follows, for $\mathbf k\in\mathbb Z^d$, we denote by $|\mathbf k|= \max_{i \in \{1, \ldots, d \}} |k_i|$.

\begin{Theorem} \label{auto-tore2}
Let $T$ be an ergodic automorphism of $\mathbb T^d$ with the notations as above.
Let $p\in ]2,4]$ and  $q$ be its conjugate exponent.
Let $f:\mathbb T^d\rightarrow \mathbb R$ be a centered function with Fourier coefficients
$(c_{\mathbf k})_{\mathbf k\in\mathbb Z^d}$ satisfying, for any integer $b \geq 2$,
\beq \label{condF1}\sum_{|\mathbf k|\ge b}|c_{\mathbf k}|^q\le R \log^{-\theta}(b) \ \text{ for some $\theta >\frac {p^2-2}{p(p-1)}$} \, ,\eeq
and
\beq \label{condF2} \sum_{|\mathbf k|\ge b}|c_{\mathbf k}|^2\le R \log^{-\beta}(b) \ \text{ for some $\beta >\frac {3p-4}{p}$} \, .\eeq
Then the series
\[\sigma^2= \bar  \lambda ((f - \bar \lambda (f))^2)
+ 2 \sum_{k>0} \bar \lambda ((f-\bar \lambda (f)) f\circ T^k)
\]
converges absolutely and, enlarging ${\mathbb T}^d$ if
necessary, there exists a sequence $(Z_i)_{i
\geq 1}$ of iid gaussian random variables with zero mean and
variance $\sigma^2$ such that, for any $t>2/p$,
\beq \label{resauto}
\sup_{1\leq k \leq n} \Big|\sum_{i=1}^kf \circ T^i - \sum_{i=1}^k Z_i\Big|  = o \big ( n^{1/p} ( \log n)^{(t+1)/2}  \big ) \text{ almost surely, as $n\rightarrow \infty$}.
\eeq
\end{Theorem}
Observe that \eqref{condF2} follows from \eqref{condF1} provided that $\theta>(3p-4)/(2p-2)$. Hence, (\ref{condF1}) and (\ref{condF2}) are both satisfied as soon as
$$ \sum_{|\mathbf k|\ge b}|c_{\mathbf k}|^q\le R \log^{-\theta}(b) \ \text{ for some $\theta >\frac {3p-4}{2(p-1)}$} \, .$$

Let us now compare our hypotheses on Fourier coefficients with those appearing
in other works.
In \cite{Leonov}, Leonov proved a central limit theorem (possibly degenerated) when
\begin{equation}\label{Leonov}
  |c_{\mathbf k}| \leq A \prod_{i=1}^d \frac{1}{(1+|k_i|)^{1/2}
  \log^{\alpha}(2+|k_i|)} \quad \text{for some $\alpha>3/2$.}
\end{equation}
In \cite{SLB}, Le Borgne proved
the functional central limit theorem and the Strassen strong invariance
principle when (\ref{condF2}) holds true with $\beta>2$ (and when $f$ is not a coboundary),
which is a weaker condition than (\ref{Leonov}).
Observe that, as $p$ converges to 2, $(p^2-2)/(p(p-1))$ and $(3p-4)/p$
both converge to 1.

\section{Probabilistic results}
\label{MR}

\setcounter{equation}{0}

In the rest of the paper, we shall use the
following notations: $\E_{k}(X)=\E (X|\mathcal{F}_{k})$, 
and  $a_{n}\ll b_{n}$ means that there exists a numerical constant $C$
not depending on $n$ such that $a_{n}\leq C b_{n}$, for all positive integers $n$.

In this section, we
give rates of convergence in the strong invariance principle under projective criteria for
stationary sequences that are non necessarily adapted to ${\mathcal
F}_i$.

\begin{Theorem} \label{ThNAS} Let $2< p < 4$ and  $t >2/p$.
Assume that $X_0$ belongs to ${\mathbb L}^p$, that
\beq \label{Cond1cob*}
\sum_{n \geq 2} \frac{n^{p-1}}{n^{2/p } (\log n)^{(t-1)p/2 }}  \big ( \Vert \E_0(X_n) \Vert^p_p  +\Vert X_{-n}- \E_0(X_{-n}) \Vert^p_p \big ) < \infty
 \, ,
\eeq
and that
\beq \label{Cond1cob**}
 \sum_{n \geq 2}   \frac{ n^{3p/4}  }{n^2  (\log n)^{(t-1)p/2 }}   \big (  \Vert \E_0(X_n  )\Vert^{p/2}_2
+  \Vert X_{-n}-\E_0(X_{-n}) \Vert^{p/2}_2 \big ) < \infty \, .
\eeq
Assume in addition that there exists a positive  integer $m$ such that
\beq \label{condcarre}
\sum_{n \geq 2} \frac{ 1 }{ n^2 (\log n)^{(t-1)p/2} }  \big \Vert  \bkE_{-n m } (S_{n}^2  ) -\bkE  (S_{n}^2)
 \big \Vert^{p/2}_{p/2} < \infty \, .\eeq
Then $n^{-1}\E(S_n^2)$ converges to  $\sigma^2=\sum_{k \in {\mathbb Z}} {\rm Cov } (X_0 , X_k)$
 and, enlarging $\Omega$ if
necessary, there exists a sequence $(Z_i)_{i
\geq 1}$ of iid Gaussian random variables with zero mean and
variance $\sigma^2$ such that
\beq \label{resThNAS}
\sup_{1\leq k \leq n} \Big|S_k - \sum_{i=1}^k Z_i\Big|  = o \big ( n^{1/p} ( \log n)^{(t+1)/2}  \big )\quad  \text{almost surely, as $n\rightarrow \infty$}.
\eeq
\end{Theorem}

\begin{Theorem} \label{ThNASp=4} Let $t >1/2$.
Assume that $X_0$ belongs to ${\mathbb L}^4$ and that the conditions (\ref{Cond1cob*}) and
(\ref{condcarre}) hold with $p=4$. Assume in addition that
\beq \label{Cond1cob**p=4}
 \sum_{n \geq 2}   n  (\log n)^{4-2t}   \big (  \Vert \E_0(X_n  )\Vert^{2}_2
+  \Vert X_{-n}-\E_0(X_{-n}) \Vert^{2}_2 \big ) < \infty \, .
\eeq
Then the conclusion of Theorem \ref{ThNAS} holds with $p=4$.
\end{Theorem}

\noindent{\bf Proof of Theorems \ref{ThNAS} and \ref{ThNASp=4}.} 
We first notice that since $p>2$, (\ref{Cond1cob*}) implies that
$$  \sum_{n>0}    n^{-1/p}   \| \E_0(X_n) \|_p  < \infty \, \text{ and } \, \sum_{n>0}    n^{-1/p}   \| X_{-n} - \E_0(X_{-n}) \|_p  < \infty $$
(apply H\"older's inequality to see this). Let $P_k(X)= \E_{k}(X)-\E_{k-1}(X)$.
Using Lemma \ref{complemmap0xi} of the appendix  with $q=1$, 
we infer that
\begin{equation}\label{Gornorm}
\sum_{k \in {\mathbb Z}} \Vert P_0(X_k) \Vert_p < \infty \, .
\end{equation}
In addition (\ref{Gornorm}) implies that $n^{-1}\E(S_n^2)$ converges to  
$\sigma^2=\sum_{k \in {\mathbf Z}} {\rm Cov } (X_0 , X_k)$.

Let now $d_0 := \sum_{j \in {\mathbb Z}} P_0 (X_j)$.
Then $d_0$ belongs to ${\mathbb L}^p$ and satisfies $\E(d_0 |{\mathcal F}_{-1})=0$.
Let $d_i := d_0 \circ T^i$ for all $i \in {\mathbb Z}$. Then $(d_i)_{i \in {\mathbb Z}}$ is a stationary sequence of martingale differences in ${\mathbb L}^p$. Let
$$ M_n := \sum_{i=1}^n d_i \, \text{ and } \, R_n := S_n -M_n \, .$$ The theorems will be proven if we can show that
\beq \label{obj1}
R_n =  o \big ( n^{1/p} ( \log n)^{(t+1)/2}  \big )\quad \text{almost surely as $n\rightarrow \infty$,}
\eeq
and that (\ref{resThNAS}) holds true with $M_k$ replacing $S_k$. Since $\E(d_0^2) = \sigma^2$ and $t>p/2$, according to
Proposition 5.1 in \cite{DDM} (applied with $\psi(n):=n^{2/p}(\log n)^t$),
to prove that (\ref{resThNAS}) holds true with $M_k$ replacing $S_k$, it suffices to prove that
\beq \label{condcarremart}
\sum_{n \geq 2} \frac{ 1 }{ n^2 (\log n)^{(t-1)p/2} }  \big \Vert  \bkE_0 (M_n^2  ) -\bkE  (M_n^2)
 \big \Vert^{p/2}_{p/2} < \infty \, .\eeq
 By standard arguments, (\ref{obj1}) will be satisfied if we can show that
 \beq \label{obj1p1}
\sum_{r >0} \frac{\Vert \max_{1 \leq \ell  \leq 2^r} |R_{\ell }| \Vert_p^p}{2^{r} \, r^{(t+1)p/2}} < \infty \, .
\eeq
Now, by stationarity,  $\Vert \max_{1 \leq \ell \leq 2^r} |R_\ell| \Vert_p \ll 2^{r/p} \sum_{k=0}^{r} 2^{-k/p} \Vert R_{2^k} \Vert_p$ (see for instance inequality (6)
in \cite{W}) and for all $i,j \geq 0$, $\| R_{i+j}\|_{q} \leq \| R_{i}\|_{q} + \| R_{j}\|_{q}$. Applying then Item 1 of Lemma 37 in \cite{MP}, we derive that
for any integer $n$ in $[2^{r},2^{r+1}[$,
\beq \label{majmaxr}
\left\Vert \max_{1 \leq \ell \leq 2^r} |R_\ell| \right\Vert_p \ll n^{1/p} \sum_{k=1}^{n} k^{-(1+1/p)} \Vert R_{k} \Vert_p \, .
\eeq
Therefore using \eqref{majmaxr} followed by an application of H\"older's inequality, we get that for any $\alpha < 1$,
\begin{align*}
\sum_{r >0} \frac{\Vert \max_{1 \leq \ell  \leq 2^r} |R_{\ell }| \Vert_p^p}{2^{r} \, r^{(t+1)p/2}} & \ll \sum_{n \geq 2} \frac{1}{n \, (\log n)^{(t+1)p/2}}
\Big ( \sum_{k=1}^{n} k^{-(1+1/p)} \Vert R_{k} \Vert_p \Big)^p \\
& \ll \sum_{n \geq 2} \frac{(\log n)^{(p-1)(1- \alpha)}}{n \, (\log n)^{(t+1)p/2}}
 \sum_{k=1}^{n} k^{-2} (\log k)^{\alpha (p-1)}\Vert R_{k} \Vert^p_p \, .
\end{align*}
Hence taking $\alpha \in ]1 -p/(2(p-1)) , 1[$ and changing the order of summation, we infer that (\ref{obj1p1}) and then (\ref{obj1}) hold provided that
\beq \label{I*}
\sum_{n \geq 1} \frac{ \Vert R_n \Vert^{p}_{p} }{ n^2 (\log n)^{(t-1)p/2} } < \infty \, .
\eeq

On an other hand, we shall prove that condition \eqref{condcarremart} is implied by: there exists  a positive finite integer $m$ such that
\beq \label{condcarremartgap}
\sum_{n \geq 2} \frac{ 1 }{ n^2 (\log n)^{(t-1)p/2} }  \big \Vert  \bkE_{-n m} (M_n^2  ) -\bkE  (M_n^2)
 \big \Vert^{p/2}_{p/2} < \infty \, .\eeq
For any nonnegative integer $i$, we set $V_i := \Vert  \bkE_0 (M_i^2  ) -\bkE  (M_i^2)
 \Vert_{p/2}$. Using that $M_n$ is a martingale, we infer that,
 for any nonnegative integers $i$ and $j$,
 \beq \label{inegaliteVi}
 V_{i+j} \leq V_i + V_j \, .
 \eeq
Let now $n \in [2^{k}, 2^{k+1}-1] \cap {\mathbb N}$, and write its binary expansion:
$$
n=\sum_{\ell=0}^{k}2^{\ell}b_{\ell}\quad \text{where }b_{k}=1\text{ and }b_{j}%
\in\{0,1\}\text{ for }j=0,\dots,k-1\,.
$$
Inequality \eqref{inegaliteVi} combined with H\"older's inequality implies that, for any $\eta >0$,
\beq \label{holder}
V_n^{p/2} \leq \Big ( \sum_{\ell=0}^{k} V_{2^\ell} \Big )^{p/2} \ll 2^{\eta p (k+1)/2}
 \sum_{\ell=0}^{k} \Big ( \frac{V_{2^\ell}}{2^{\eta \ell } }\Big )^{p/2} \, .\eeq
Therefore
$$
\sum_{n \geq 2} \frac{ 1 }{ n^2 (\log n)^{(t-1)p/2} }  V_n^{p/2} \ll \sum_{k >0}
  \frac{ 2^{\eta p (k+1)/2} }{2^k k^{(t-1)p/2} }   \sum_{\ell=0}^{k}
   \Big ( \frac{V_{2^\ell}}{2^{\eta \ell } }\Big )^{p/2}\, .
$$
Changing the order of summation and taking $\eta \in ]0,2/p[$,
it follows that \eqref{condcarremart} is implied by \beq \label{condcarremart*}
\sum_{k \geq 1} \frac{ 1 }{ 2^k k^{(t-1)p/2} }  \big \Vert  \bkE_0 (M_{2^k}^2  ) -\bkE  (M_{2^k}^2)
 \big \Vert^{p/2}_{p/2} < \infty \eeq
(actually due to the subadditivity of the sequence $(V_i)$ both conditions are equivalent, see the proof of item 1 of Lemma 37 in \cite{MP} to prove that \eqref{condcarremart} entails \eqref{condcarremart*}).
Now, since $(M_n)$ is a martingale,
 $$
  \bkE_0 (M_{2^k}^2  ) -\bkE  (M_{2^k}^2) = \sum_{j=1}^k \big(  \bkE_0 ((M_{2^j}-M_{2^{j-1}})^2  ) -\bkE  ((M_{2^j}-M_{2^{j-1}})^2) \big ) + \bkE_0 (d_1^2) - \E(d_1^2) \, ,
 $$
which implies by stationarity that
 $$
 \big \Vert  \bkE_0 (M_{2^k}^2  ) -\bkE  (M_{2^k}^2)
 \big \Vert_{p/2} \leq \sum_{j=0}^{k-1}  \big \Vert  \bkE_{-2^j} (M_{2^j}^2  ) -\bkE  (M_{2^j}^2)
 \big \Vert_{p/2} +  \big \Vert   \bkE_0 (d_1^2) - \E(d_1^2)  \big \Vert_{p/2} \, .
 $$
 Therefore by using H\"older's inequality as done in \eqref{holder} with $\eta \in ]0,2/p[$, we infer that \eqref{condcarremart*} is implied by
 \beq \label{condcarremart**}
\sum_{k \geq 1} \frac{ 1 }{ 2^k k^{(t-1)p/2} }  \big \Vert  \bkE_{-2^k} (M_{2^k}^2  ) -\bkE  (M_{2^k}^2)
 \big \Vert^{p/2}_{p/2} < \infty \, .\eeq
 Notice now that the sequence $(W_n)_{n>0}$ defined  by
$$W_n := \big \Vert  \bkE_{-n} (M_{n}^2  ) -\bkE  (M_{n}^2)
 \big \Vert_{p/2} $$ is subadditive. Indeed,  for any non negative integers $i$ and $j$,
using that $M_n$ is a martingale together with the stationarity, we derive that
  \begin{align*}
  W_{i+j} & = \big \Vert  \bkE_{-(i+j)} (M_{i}^2  ) -\bkE  (M_{i}^2) +  \bkE_{-(i+j)} ((M_{i+j} -M_{i})^2  ) -\bkE  ((M_{i+j} -M_{i})^2)
 \big \Vert_{p/2} \\
 & \leq \big \Vert  \bkE_{-i} (M_{i}^2  ) -\bkE  (M_{i}^2)\big \Vert_{p/2} +  \big \Vert \bkE_{-j} (M_{j} )^2  ) -\bkE  (M_{j})^2)
 \big \Vert_{p/2} \\
 & \leq W_i + W_j\, .
  \end{align*}
Therefore $W_{i+j}^{p/2} \leq 2^{p/2} W_i^{p/2} + 2^{p/2} W_j^{p/2}$. This implies that,
for any integer $\ell$ and any
integer $0\leq j\leq \ell$,
\begin{equation}
W_{\ell}^{p/2}\leq 2^{p/2} ( W^{p/2}_{j}+W^{p/2}_{\ell-j} ), \ \text{in such a way that }(\ell +1)W^{p/2}_{\ell}\leq 2^{1+ p/2} \sum_{j=1}%
^{\ell}W^{p/2}_{j}\,. \label{condsubadd2}%
\end{equation}
Therefore using the second part of \eqref{condsubadd2} with $\ell=2^k$,
we infer that condition \eqref{condcarremart**} is implied by
  \beq \label{condcarremart***}
\sum_{n \geq 2} \frac{ 1 }{n^2 (\log n)^{(t-1)p/2} }  \big \Vert  \bkE_{-n} (M_{n}^2  ) -\bkE  (M_{n}^2)
 \big \Vert^{p/2}_{p/2} < \infty \, .\eeq
It remains to prove that \eqref{condcarremartgap} implies \eqref{condcarremart***}.
With this aim, we have, for any positive integer $m$,
$$
M_n = \sum_{k=1}^{m} \big ( M_{k[n m^{-1}]} - M_{(k-1)[n m^{-1}]}\big ) + M_n -M_{m[n m^{-1}]} \, .
$$
Using that $M_n$ is a martingale together with the stationarity, we then infer that
\begin{multline*}
 \big \Vert  \bkE_{-n} (M_{n}^2  ) -\bkE  (M_{n}^2)
 \big \Vert^{p/2}_{p/2} \leq 2^{p/2} m^{p/2}  \big \Vert  \bkE_{-n} (M_{[n m^{-1}]}^2  ) -\bkE  (M_{[n m^{-1}]}^2)
 \big \Vert^{p/2}_{p/2} \\
 + 2^{p/2} \big \Vert  \bkE_{-n} (M_{n-m[n m^{-1}]}^2  ) -\bkE  (M_{n-m[n m^{-1}]}^2)
 \big \Vert^{p/2}_{p/2} \, ,
\end{multline*}
which, together with the fact that $ n-m[n m^{-1}] < m  $, implies that
\begin{align} \label{decgap1}
 \big \Vert  \bkE_{-n} (M_{n}^2  ) & -\bkE  (M_{n}^2)
 \big \Vert^{p/2}_{p/2} \leq 2^{p/2} m^{p/2} \Big ( 2^{p/2}\Vert d_0 \Vert_p^{p} + \big \Vert  \bkE_{-n} (M_{[n m^{-1}]}^2  ) -\bkE  (M_{[n m^{-1}]}^2)
 \big \Vert^{p/2}_{p/2} \Big) \nonumber \\
 & \leq 2^{p/2} m^{p/2} \Big ( 2^{p/2}\Vert d_0 \Vert_p^{p} + \big \Vert  \bkE_{-m[n m^{-1}]} (M_{[n m^{-1}]}^2  ) -\bkE  (M_{[n m^{-1}]}^2)
 \big \Vert^{p/2}_{p/2} \Big)\, ,
\end{align}
where for the last line we have used the fact that $n \geq m[n m^{-1}]$.
We notice now that due to the martingale property of $(M_n)$ and to stationarity, the sequence $ (U_i)_{i \geq 0}$ defined for any non negative integer $i$ by
$$
U_i := \big \Vert  \bkE_{-mi} (M_{i}^2  ) -\bkE  (M_{i}^2)
 \big \Vert^{p/2}_{p/2}
$$
satisfies, for any positive integers $i$ and $j$,
\begin{multline*}
U_{i+j}  \leq \Big ( \big \Vert  \bkE_{-m(i+j)} (M_{i}^2  ) -\bkE  (M_{i}^2)
 \big \Vert_{p/2} \\
 + \big \Vert  \bkE_{-m(i+j)} ((M_{i+j} -M_i)^2  ) -\bkE  ((M_{i+j} -M_i)^2 )
 \big \Vert_{p/2}  \Big )^{p/2}
  \leq 2^{p/2} U_i + 2^{p/2} U_j \, .
\end{multline*}
Hence by (\ref{condsubadd2}) applied with $W_i^{p/2} = U_{i}$,
\beq \label{decgap2} U_{[n m^{-1}] } \leq 2^{1+ p/2} ([n m^{-1}] + 1)^{-1} \sum_{k=1}^{[n m^{-1}]} U_k \leq 2^{1+ p/2}  \sum_{k=1}^{[n m^{-1}]} \frac{U_k}{k} \, .
\eeq
Therefore starting from (\ref{decgap1}), considering (\ref{decgap2}) and changing the order of summation,
we infer that \eqref{condcarremart***} (and so \eqref{condcarremart}) holds provided that  \eqref{condcarremartgap} does. To end the proof, it remains to show that under the conditions of Theorems \ref{ThNAS} and \ref{ThNASp=4}, the conditions (\ref{I*}) and \eqref{condcarremartgap} are satisfied. This is achieved by using the two following lemmas.

\begin{Lemma} \label{normrn} Let $p \in [2,4]$. Assume that (\ref{Cond1cob*}) holds. Then
$$
\sum_{n \geq 1} \frac{ \max_{1 \leq \ell \leq n}\Vert R_{\ell} \Vert^{p}_{p} }{ n^2 (\log n)^{(t-1)p/2} } < \infty \, ,
$$
and (\ref{I*}) holds.
\end{Lemma}
\begin{Lemma} \label{normMn} Let $p \in [2,4]$ and assume that (\ref{Cond1cob*})  and (\ref{condcarre}) are satisfied. Assume in addition that (\ref{Cond1cob**}) holds when $2 < p< 4 $ and (\ref{Cond1cob**p=4}) does when $p=4$ . Then \eqref{condcarremartgap} is satisfied.
\end{Lemma}
It remains to prove the two above lemmas.

\medskip

\noindent{\bf Proof of Lemma \ref{normrn}.} Since (\ref{Cond1cob*}) implies (\ref{Gornorm}),  Item 2 of Proposition \ref{approxrnq}
given in the appendix implies that, for any positive integers ${\ell}$ and  $N$,
$$
\Vert R_{\ell} \Vert_{p} \ll   \max_{k={\ell},N}\Vert \E_0 (S_{k}  ) \Vert_p + \max_{k={\ell},N}    \Vert S_k - {\mathbb E}_k (S_k )  \Vert_p
 + {\ell}^{1/2}  \sum_{|j | \geq N} \Vert    P_0(X_j) \Vert_p
 \, .
$$
Next, applying Lemma \ref{complemmap0xi} given in the appendix with $q=1$, and using the fact that by stationarity, for any positive integer $k$,
\beq \label{lma1p0}
\Vert  \E_0 (S_{k} ) \Vert_p \leq \sum_{\ell =1}^{k} \Vert \E_0(X_{\ell}) \Vert_p \, \text{ and } \,  \Vert S_{k} - \E_{k} (S_{k} ) \Vert_p \leq \sum_{\ell=0}^{k-1}\Vert X_{-\ell} - \E_0 (X_{-\ell} ) \Vert_p \, , \eeq
we derive that for, any positive integers $N \geq n$,
\begin{multline} \label{lma1p1}
\max_{1 \leq \ell \leq n}\Vert R_{\ell} \Vert_{p} \ll    \sum_{k=1}^{N} \Vert \E_0(X_k) \Vert_p +
 \sum_{k=0}^{N-1}\Vert X_{-k} - \E_0 (X_{-k} ) \Vert_p +\\
+ n^{1/2}  \sum_{k \geq [N/2]  } \frac{\Vert  \E_0(X_{k} )\Vert_p}{k^{1/p}} + n^{1/2} \sum_{k \geq [N/2]  } \frac{\Vert X_{-k} - \E_0(X_{-k} )\Vert_p}{k^{1/p}}
 \, .
\end{multline}
The lemma follows from (\ref{lma1p1}) with $N=[n^{p/2}]$ by using H\"older's
inequality (see the computations  in the proof of Proposition 2.2 in \cite{DDM}). $\square$

\medskip

\noindent{\bf Proof of Lemma \ref{normMn}.}
Let $m$ be a positive integer such that (\ref{condcarre}) is satisfied. We first write that
\begin{multline*}
\Vert  \bkE_{ -nm } (M_n^2  ) -\bkE  (M_n^2)
 \big \Vert_{p/2} \leq \Vert  \bkE_{ -nm } (S_n^2  ) -\bkE  (S_n^2)
 \big \Vert_{p/2} \\
 +2 \Vert \bkE_{- nm }  (S_nR_n  ) -
\bkE  (S_nR_n)  \Vert_{p/2} + 2 \Vert R_n \Vert^{2}_{p} \, .
\end{multline*}
By using Lemma \ref{normrn}, and since \eqref{condcarre} holds, Lemma \ref{normMn} will follow if we can prove that
\beq \label{II*}
\sum_{n \geq 1} \frac{ 1 }{  n^2 (\log n)^{(t-1)p/2} } \Vert \E_{ -nm }  (S_nR_n  )\Vert^{p/2}_{p/2} < \infty\, .\eeq
With this aim we shall prove the following inequality. For any non negative integer $r$ and any positive integer $u_n$ such that $u_n \leq n$, we have that
\begin{multline} \label{doubleprod}
\Vert \E_{-r}  (S_nR_n  )\Vert_{p/2} \ll \sqrt{u_n} \big ( \Vert \E_{0}  (S_n  )\Vert_{2} +
 \Vert S_n - \E_n  (S_n  )\Vert_{2} \big ) + \max_{k=\{n,n-u_n\}} \Vert R_k \Vert_p^2 +\\
  +\sqrt{n} \big ( \Vert \E_{-u_n}  (S_n  )\Vert_{2} + \Vert S_n - \E_{n+u_n}  (S_n  )\Vert_{2} \big )  + \max_{k=\{n,u_n\}}\Vert \E_{-r}  (S^2_k  ) - \E(S_k^2)\Vert_{p/2} \\
 + \sqrt{n} \Big (  \sum_{k=1}^n
\big \Vert  \sum_{|j| \geq k+n}  P_0(X_j) \big \Vert_2^{2} \Big )^{1/2} \, .
\end{multline}
Let us show how, thanks to (\ref{doubleprod}), the convergence (\ref{II*}) can be proven. Let us first consider the case where $2<p<4$. Notice that the following
elementary claim is valid:
\begin{Claim} \label{claimfacile}
If ${\mathcal F}$ and
${\mathcal G}$ are two $\sigma$-algebras such that ${\mathcal G }\subset {\mathcal F}$, then
for any random variable $X$ in ${\mathbb L}^q$ for  $q \geq 1$, $\Vert X - {\mathbb E}(X | {\mathcal F}) \Vert_q \leq
2 \Vert X - {\mathbb E}(X | {\mathcal G}) \Vert_q$.
\end{Claim}
Starting from (\ref{doubleprod}) with $r=nm$ and $u_n=n$, and using Claim \ref{claimfacile}, we derive that
\begin{multline*}
\Vert \E_{-nm}  (S_nR_n  )\Vert_{p/2} \ll  \Vert \E_{-nm} (S^2_n  ) - \E(S_n^2)\Vert_{p/2}
+ \sqrt{n} \big ( \Vert \E_0  (S_n  )\Vert_{2} + \Vert S_n - \E_n  (S_n  )\Vert_{2} \big )+ \\
+ \Vert R_n \Vert_p^2 + n
 \sum_{|j| \geq n} \Vert  P_0(X_j) \Vert_2  \, .
\end{multline*}
This last inequality combined with condition (\ref{condcarre}) and Lemma \ref{normrn} shows that (\ref{II*}) will be satisfied if we can prove that
\beq \label{II*1}
\sum_{n \geq 1} \frac{ n^{p/4} }{  n^2 (\log n)^{(t-1)p/2} } \big ( \Vert \E_0  (S_n  )\Vert_{2} + \Vert S_n - \E_n  (S_n  )\Vert_{2} \big )^{p/2} < \infty\, ,\eeq
and
\beq \label{II*2}
\sum_{n \geq 1} \frac{ n^{p/2} }{  n^2 (\log n)^{(t-1)p/2} } \big ( \sum_{|j| \geq n} \Vert  P_0(X_j) \Vert_2 \big )^{p/2} < \infty\, .\eeq
To prove \eqref{II*1}, we use the inequalities (\ref{lma1p0}) with $p=2$. Hence setting
\beq \label{defal}
a_{\ell} = \Vert \E_0(X_{\ell}) \Vert_2 + \Vert X_{-\ell+1} - \E_0 (X_{-\ell+1} ) \Vert_2 \, ,
\eeq
and using H\"older's inequality, we derive that for any $\alpha <1$,
\begin{multline*}
\sum_{n \geq 1} \frac{ n^{p/4} }{  n^2 (\log n)^{(t-1)p/2} } \big ( \Vert \E_0  (S_n  )\Vert_{2} + \Vert S_n - \E_n  (S_n  )\Vert_{2} \big )^{p/2}
\ll \sum_{n \geq 1} \frac{ n^{p/4} }{  n^2 (\log n)^{(t-1)p/2} } \Big ( \sum_{\ell=1}^n a_{\ell} \Big )^{p/2}\\
\ll \sum_{n \geq 1} \frac{ n^{p/4} n^{(1-\alpha)(p/2-1)}}{  n^2 (\log n)^{(t-1)p/2} } \sum_{\ell=1}^n \ell^{\alpha (p/2-1)}a_{\ell}^{p/2} \, .
\end{multline*}
Taking $\alpha \in ](3p-8)/(2p-4), 1[$ (this is possible since $p <4$) and changing the order of summation, we infer that \eqref{II*1} holds provided that
\eqref{Cond1cob**} does. It remains to show that \eqref{II*2} is satisfied. Using Lemma \ref{complemmap0xi} and the notation \eqref{defal}, we first observe that
$$
\sum_{n \geq 1} \frac{ n^{p/2} }{  n^2 (\log n)^{(t-1)p/2} } \big ( \sum_{|j| \geq n} \Vert  P_0(X_j) \Vert_2 \big )^{p/2} \ll \sum_{n \geq 1} \frac{ n^{p/2} }{  n^2 (\log n)^{(t-1)p/2} } \Big ( \sum_{\ell \geq [n/2]} \ell^{-1/2}a_{\ell}\Big)^{p/2} \, .
$$
Therefore by H\"older's inequality, it follows that for any $\alpha >1$,
$$
\sum_{n \geq 1} \frac{ n^{p/2} }{  n^2 (\log n)^{(t-1)p/2} } \big ( \sum_{|j| \geq n} \Vert  P_0(X_j) \Vert_2 \big )^{p/2} \ll \sum_{n \geq 1} \frac{ n^{p/2} n^{(1-\alpha)(p/2 -1)}}{  n^2 (\log n)^{(t-1)p/2} } \sum_{\ell \geq [n/2]} \ell^{\alpha(p/2 -1)}\ell^{-p/4}a_{\ell}^{p/2} \, .
$$
Therefore taking $\alpha \in ]1,2[$ and changing the order of summation, we infer that \eqref{II*1} holds provided that
\eqref{Cond1cob**} does. This ends the proof of  (\ref{II*}) when $p \in ]2,4[$.

 Now, we prove (\ref{II*}) when $p=4$. With this aim we start from (\ref{doubleprod}) with $r= nm$ and $u_n=[\sqrt n]$. This  inequality combined with condition (\ref{condcarre}), Lemma \ref{normrn} and the arguments developed to prove \eqref{II*1} and \eqref{II*2} shows that (\ref{II*}) will be satisfied for $p=4$ if we can prove that
\beq \label{II*11}
\sum_{n \geq 1} \frac{ 1}{  n(\log n)^{2(t-1)} } \big ( \Vert \E_{-[\sqrt n ]}  (S_n  )\Vert_{2} + \Vert S_n - \E_{n +[\sqrt n ]} (S_n  )\Vert_{2} \big )^2 < \infty\, ,\eeq
and
\beq \label{II*22}
\sum_{n \geq 2} \frac{ 1 }{ n^2 (\log n)^{2(t-1)} }  \big \Vert  \E_{-nm} (S_{[\sqrt n]}^2  ) -\E  (S_{[\sqrt n]}^2)
 \big \Vert^{2}_{2} < \infty \, . \eeq
We start by proving \eqref{II*11}. With this aim, using the notation \eqref{defal}, we first write that
$$
\Vert \E_{-[\sqrt n ]}  (S_n  )\Vert_{2} + \Vert S_n - \E_{n+[\sqrt n ]}  (S_n  )\Vert_{2}\leq \sum_{k=[\sqrt n ]+1}^{n+[\sqrt n ]} a_k \, .
$$
Therefore by Cauchy-Schwarz's inequality
\begin{multline*}
\sum_{n \geq 1} \frac{ 1}{  n(\log n)^{2(t-1)} } \big ( \Vert \E_{-[\sqrt n ]}  (S_n  )\Vert_{2} + \Vert S_n - \E_{n +[\sqrt n ]} (S_n  )\Vert_{2} \big )^2  \ll
\sum_{n \geq 1} \frac{ \log n }{  n(\log n)^{2(t-1)} } \sum_{k=[\sqrt n ]+1}^{n+[\sqrt n ]} k a_k^2 \\
\ll \sum_{n \geq 1} \frac{ 1}{  n } \sum_{k=[\sqrt n ]+1}^{n+[\sqrt n ]} \frac{k\, \log k}{(\log k)^{2(t-1)}} a_k^2 \, .
\end{multline*}
Changing the order of summation, this proves that \eqref{II*11} holds provided that (\ref{Cond1cob**p=4}) does. It remains to prove \eqref{II*22}. With this aim, we set for any positive real $x$,
$$
h ([x]) = \big \Vert  \E_{-m[x]} (S_{[x]}^2  ) -\E  (S_{[x]}^2)
 \big \Vert^{2}_{2} \, ,
$$
and we notice that, for any integer $n\geq 0$, $\big \Vert  \E_{-nm} (S_{[\sqrt n]}^2  ) -\E  (S_{[\sqrt n]}^2)
 \big \Vert^{2}_{2} \leq h([\sqrt n])$. In addition, if $x \in [n, n+1[$, then $[\sqrt n]=[\sqrt x]$ or $[\sqrt{ n}]=[\sqrt x]-1$. Therefore
 \begin{align*}
\sum_{n \geq 3} & \frac{ 1 }{ n^2 (\log n)^{(t-1)p/2} } h([\sqrt n]) \ll  \sum_{n \geq 3} h([\sqrt n]) \int_{[n, n+1[} \frac{ 1 }{ x^2 (\log x)^{(t-1)p/2} }  dx \\
&\ll \int_{3}^{\infty} \frac{ 1 }{ x^2 (\log x)^{(t-1)p/2} } h([\sqrt x] ) dx + \int_{3}^{\infty} \frac{ 1 }{ x^2 (\log x)^{(t-1)p/2} } h([\sqrt x] -1) dx \\
&\ll \int_{2}^{\infty} \frac{ 1 }{ y^3 (\log y)^{(t-1)p/2} } h([y] ) dy \ll \sum_{n \geq 2} \frac{ 1 }{ n^3 (\log n)^{(t-1)p/2} } h(n ) dy  \, .
 \end{align*}
For the last inequality, we have used that if $y \in [n, n+1[$, then $[y] =n$.
Therefore condition (\ref{condcarre}) implies \eqref{II*22}. This ends the proof of \eqref{II*} when $p=4$.

\medskip
It remains to prove (\ref{doubleprod}). With this aim, we start with the decomposition of $R_n$ given in Proposition \ref{approxrnq}
of the appendix with $N=n$. Therefore setting
$$
A_n := \sum_{k=1}^n \sum_{j \geq 2n +1} P_k(X_j) +\sum_{k=1}^n \sum_{j\geq n} P_k(X_{-j}) \, ,
$$we write that
\beq \label{p1double} R_n   =   \E_0 (S_n ) -\E_0 (S_{n} ) \circ T^n + \E_{-n} (S_{n}) \circ T^n
 + S_n - {\mathbb E}_n (S_n ) - (  \E_{2n} (  S_{n} - \E_n (S_{n} ) ) \circ T^{-n}  -A_n  \, .
\eeq
Starting from \eqref{p1double} and noticing that $$ \Vert \E_{-r}(S_n (\E_{-n} (S_{n}) \circ T^n  ) \Vert_{p/2} \leq  \Vert \E_{0}(S_n (\E_{-n} (S_{n}) \circ T^n  ) \Vert_{p/2} \leq \Vert \E_0 (S_n) \Vert_{p} \Vert \E_0 (S_{2n} -S_n) \Vert_p \, ,
$$
and that $\E_{-r}(S_n (S_n - \E_n (S_n))= \E_{-r} ( (S_n - \E_n(S_n))^2)$, we first get
\begin{multline} \label{supstep}
 \Vert \E_{-r}  (S_nR_n ) \Vert_{p/2} \leq 2 \Vert \E_0 (S_n) \Vert_{p}^2 + \Vert S_n-
\E_n (S_n) \Vert_{p}^2 + \Vert \E_{-r} (S_n\E_n (S_{2n}-S_n)) \Vert_{p/2}  \\
 + \Vert \E_{-r} (S_n\E_n (S_n \circ T^{-n} - \E_{0} (S_n \circ T^{-n})) )\Vert_{p/2} + \Vert  \E_{-r} (S_{n}A_n)  \big \Vert_{p/2} \, .
\end{multline}
Next, we use the following fact: if $X$ and $Y$ are two variables in ${\mathbb L}^p$ with $p \in [2,4]$, then for any integer  $u$,
\beq \label{fact}
\Vert \E_u(XY) \Vert_{p/2} \leq \Vert \E_u(X^2) - \E(X^2) \Vert_{p/2} + \Vert Y \Vert_{p}^2 + \sqrt{\E(X^2)} \Vert Y \Vert_{2} \, .
\eeq
Indeed, it suffices to write that
\begin{align*}
 \Vert \E_u  (XY ) \Vert_{p/2} & \leq  \Vert \E^{1/2}_u  (X^2) \E^{1/2}_u(Y ^2  ) \Vert_{p/2} \nonumber \\
 & \leq   \Vert | \E_u  (X^2)- \E(X^2) |^{1/2}\E^{1/2}_u(Y^2) \Vert_{p/2}+  (\E(X^2) )^{1/2}\Vert \E^{1/2}_u(Y^2 ) \Vert_{p/2}  \nonumber \\
 & \leq   \Vert \E_u  (X^2)- \E(X^2)  \Vert_{p/2} + \Vert Y \Vert^2_{p} +  (\E(X^2) )^{1/2}\Vert \E^{1/2}_u(Y^2 ) \Vert_{p/2} \, ,
\end{align*}
and to notice that, since $p \in [2,4]$, 
$\Vert \E^{1/2}_u(Y^2 ) \Vert_{p/2}\leq \Vert \E^{1/2}_u(Y^2 ) \Vert_{2}=\Vert Y \Vert_{2}$.
Therefore, starting from \eqref{supstep} and using (\ref{fact}) together with
 $\E(S_n^2) \ll n$,  we infer that
\begin{multline*}
 \Vert \E_{-r}  (S_nR_n ) \Vert_{p/2} \ll \Vert \E_0 (S_n) \Vert_{p}^2 + \Vert S_n- \E_n (S_n)
 \Vert_{p}^2 + \Vert \E_{-r} (S_n\E_n (S_{2n}-S_n)) \Vert_{p/2} \\
 + \Vert \E_{-r} (S_n\E_n (S_n \circ T^{-n} - \E_{0} (S_n \circ T^{-n}))) \Vert_{p/2}\\ +
\Vert  \E_{-r} (S_{n}^2  ) -\E  (S_{n}^2)  \big \Vert_{p/2}
 + \Vert A_n \Vert_p^2 + n^{1/2} \Vert A_n \Vert_2 \, ,
\end{multline*}
and since $\Vert \E_0 (S_n) \Vert_{p} \leq \Vert R_n \Vert_p$, $ \Vert S_n- \E_n (S_n) \Vert_{p} \leq 2 \Vert R_n \Vert_p $ and $\Vert A_n \Vert_p \leq 8 \Vert R_n \Vert_p$, we have overall that
\begin{multline} \label{p2double}
 \Vert \E_{-r}  (S_nR_n ) \Vert_{p/2} \ll \Vert R_n \Vert_p^2 +  \Vert \E_{-r}
(S_n\E_n (S_{2n}-S_n)) \Vert_{p/2} \\
 + \Vert \E_{-r} (S_n\E_n (S_n \circ T^{-n} - \E_0 (S_n \circ T^{-n})) )\Vert_{p/2} + \Vert
 \E_{-r} (S_{n}^2  ) -\E  (S_{n}^2)  \big \Vert_{p/2} +  n^{1/2} \Vert A_n \Vert_2 \, .
\end{multline}
By orthogonality and by stationarity,
\begin{eqnarray}
\Vert A_n \Vert_2  &\leq&  \Big (  \sum_{k=1}^n \Big \Vert \sum_{j \geq 2n +1} P_k(X_j)
   \Big \Vert^2_2  \Big )^{1/2} +
     \Big ( \sum_{k=1}^n \Big \Vert \sum_{j\geq n} P_k(X_{-j})
    \Big \Vert^2_2  \Big )^{1/2}\nonumber\\
&\leq&  \Big(  \sum_{k=1}^n \Big \Vert \sum_{\ell \geq k+n} P_0(X_\ell) \Big \Vert^2_2  \Big )
^{1/2} + \Big ( \sum_{k=1}^n \Big \Vert \sum_{\ell \geq k +n} P_0(X_{-\ell}) \Big
   \Vert^2_2  \Big )^{1/2}  \, .\label{p2doublebis}
\end{eqnarray}
Now for any integer $u_n$  such that $u_n \leq n$,
\begin{multline} \label{p3double}
\Vert \E_{-r} (S_n\E_n (S_{2n}-S_n)) \Vert_{p/2}\leq \Vert \E_{-r} ((S_{n} - S_{n-u_n})\E_n (S_{2n}-S_n)) \Vert_{p/2}\\ + \Vert \E_{-r} (S_{n-u_n}\E_n (S_{2n}-S_n)) \Vert_{p/2}  \\
\ll \Vert  \E_{-r} (S_{u_n}^2  ) -\E  (S_{u_n}^2)  \big \Vert_{p/2} +  \Vert \E_{0}(S_n) \Vert_p^2 + \sqrt{u_n}\Vert \E_0(S_n) \Vert_2\\+\Vert \E_{-r} (S_{n-u_n}\E_n (S_{2n}-S_n)) \Vert_{p/2}  \, ,
\end{multline}
where for the last inequality we have used (\ref{fact}) together with $\E(S_{u_n}^2) \ll u_n$.
Next, we write that
\begin{multline*}
 \Vert \E_{-r} (S_{n-u_n}\E_n (S_{2n}-S_n)) \Vert_{p/2}
 \leq \Vert \E_{-r} ((S_{n-u_n} - \E_{n-u_n} (S_{n-u_n}))\E_n (S_{2n}-S_n)) \Vert_{p/2} \\
+\Vert \E_{-r} (\E_{n-u_n} (S_{n-u_n})\E_n (S_{2n}-S_n)) \Vert_{p/2}  \\
 \leq \Vert S_{n-u_n} - \E_{n-u_n} (S_{n-u_n})\Vert_{p}  \Vert \E_0(S_n) \Vert_p +\Vert \E_{-r} (\E_{n-u_n} (S_{n-u_n})\E_{n-u_n} (S_{2n}-S_n)) \Vert_{p/2} \\
 \leq \Vert S_{n-u_n} - \E_{n-u_n} (S_{n-u_n})\Vert^2_{p} +  \Vert \E_0(S_n) \Vert^2_p
 +\Vert \E_{-r} (S_n\E_{n-u_n} (S_{2n}-S_n)) \Vert_{p/2}  \\ + \Vert \E_{-r} ((S_n - S_{n-u_n})\E_{n-u_n} (S_{2n}-S_n)) \Vert_{p/2} \, .
\end{multline*}
Therefore using \eqref{fact}, we infer that
\begin{multline} \label{p4double}
\Vert \E_{-r} (S_{n-u_n}\E_n (S_{2n}-S_n)) \Vert_{p/2} \ll  \max_{k=\{n,n-u_n\}} \Vert R_k \Vert_p^2 + \sqrt{n} \Vert \E_{-u_n} (S_n) \Vert_2 \\
+ \max_{k=\{n,u_n\}}\Vert \E_{-r}  (S^2_k  ) - \E(S_k^2)\Vert_{p/2}
 \, .
\end{multline}
We deal now with the third term in the right-hand side of (\ref{p2double}). With this aim, we first write that
\begin{multline}  \label{p5double}
\Vert \E_{-r} (S_n\E_n (S_n \circ T^{-n} - \E_0 (S_n \circ T^{-n}))) \Vert_{p/2} \leq \Vert \E_{-r} (S_n\E_n (S_n \circ T^{-n} - \E_{u_n}(S_n \circ T^{-n})) )\Vert_{p/2} \\ + \Vert \E_{-r} (S_n\E_{u_n} (S_n \circ T^{-n} - \E_{0}(S_n \circ T^{-n})) )\Vert_{p/2}.
\end{multline}
By using \eqref{fact} together with $\E (S_{u_n}^2 )\ll n$, stationarity and the fact that $\Vert S_n - \E_{n+u_n} (S_n) \Vert_2 \leq 2 \Vert R_n \Vert_p$, we infer that
\begin{multline}  \label{p6double}
\Vert \E_{-r} (S_n\E_n (S_n \circ T^{-n} - \E_{u_n}(S_n \circ T^{-n})) )\Vert_{p/2} \ll \Vert \E_{-r}  (S^2_n  ) - \E(S_n^2)\Vert_{p/2} \\
+ \Vert R_n \Vert_p^2 + \sqrt{n} \Vert S_n - \E_{n+u_n} (S_n) \Vert_2 \, .
\end{multline}
On the other hand,
\begin{multline*}
\Vert \E_{-r} (S_n\E_{u_n} (S_n \circ T^{-n} - \E_{0}(S_n \circ T^{-n})) ) \Vert_{p/2} \leq \Vert \E_{-r} (S_{u_n}\E_{u_n} (S_n \circ T^{-n} - \E_{0}(S_n \circ T^{-n}))) \Vert_{p/2}  \\
+ \Vert \E_{-r} (\E_{u_n} (S_n -S_{u_n})\E_{u_n} (S_n \circ T^{-n} - \E_{0}(S_n \circ T^{-n})) \Vert_{p/2}\, .
\end{multline*}
We apply (\ref{fact}) to the first term of the right hand side together with $\E (S_{u_n}^2) \ll n$. Hence by stationarity and since
$\Vert S_n - \E_{n}(S_n ) \Vert_{p} \leq 2 \Vert R_n \Vert_{p}$, we derive that
\begin{multline*}
 \Vert \E_{-r} (S_{u_n}\E_{u_n} (S_n \circ T^{-n} - \E_{0}(S_n \circ T^{-n})) ) \Vert_{p/2} \ll
\Vert \E_{-r}  (S^2_{u_n}  ) - \E(S_{u_n}^2)\Vert_{p/2}+ \\+ \Vert R_n \Vert_p^2 + \sqrt{u_n} \Vert S_n - \E_{n} (S_n) \Vert_2 \, .
\end{multline*}
On the other hand, by stationarity,
\begin{align*}
\Vert \E_{-r} (\E_{u_n} (S_n -S_{u_n}) & \E_{u_n} (S_n \circ T^{-n} - \E_{0}(S_n \circ T^{-n})) )
 \Vert_{p/2}  \\
& \leq \Vert \E_{u_n} (S_n -S_{u_n}) \Vert_{p} \Vert \E_{u_n} (S_n \circ T^{-n} - \E_{0}(S_n \circ T^{-n})) \Vert_{p}
\\
& \leq \Vert \E_{0} (S_{n -u_n}) \Vert_{p}\Vert S_n - \E_{n}(S_n )) \Vert_{p} \, .
\\
& \leq \Vert \E_{0} (S_{n -u_n}) \Vert_{p}^2 +  \Vert R_n \Vert^2_{p} \, .
\end{align*}
Therefore we get overall that
\begin{multline}  \label{p7double}
\Vert \E_{-r} (S_n\E_{u_n} (S_n \circ T^{-n} - \E_{0}(S_n \circ T^{-n})) )\Vert_{p/2} \ll \Vert R_n \Vert_p^2 + \Vert \E_0(S_{n-u_n}) \Vert_p^2 \\
+ \Vert \E_{-r}  (S^2_{u_n}  ) - \E(S_{u_n}^2)\Vert_{p/2}  +  \sqrt{u_n} \Vert S_n - \E_{n} (S_n) \Vert_2
 \, .
\end{multline}
Starting from (\ref{p5double}) and taking into account (\ref{p6double}) and (\ref{p7double}), we get that
\begin{multline}  \label{p8double}
\Vert \E_{-r} (S_n\E_n (S_n \circ T^{-n} - \E_0 (S_n \circ T^{-n})) )\Vert_{p/2} \ll \sqrt{u_n}  \Vert S_n - \E_n  (S_n  )\Vert_{2}  +\sqrt{n} \Vert S_n - \E_{n+u_n}  (S_n  )\Vert_{2}  \\ + \max_{k=\{n,u_n\}}\Vert \E_{-r}  (S^2_k  ) - \E(S_k^2)\Vert_{p/2}
+ \max_{k=\{n,n-u_n\}} \Vert R_k \Vert_p^2 \, .
\end{multline}

Finally, starting from \eqref{p2double} and considering \eqref{p2doublebis}, \eqref{p3double}, \eqref{p4double} and \eqref{p8double},
we conclude that (\ref{doubleprod}) holds.
 $\square$

\section{Proof of Theorem  \ref{auto-tore2}}

\subsection{Preparatory material}

Let us denote by $E_u$, $E_e$ and $E_s$ the $S$-stable vector spaces associated to the eigenvalues
of $S$ of modulus respectively larger than one, equal to one and smaller than one.
Let $d_u$, $d_e$ and $d_s$ be their respective dimensions.
Let $v_1,...,v_d$ be a basis of $\mathbb R^d$ in which
$S$ is represented by a real Jordan matrix.
Suppose that $v_1,...,v_{d_u}$ are in
$E_u$, $v_{d_u+1},...,v_{d_u+d_e}$ are in $E_e$ and $v_{d_u+d_e+1},...,v_d$ are in $E_s$.
We suppose moreover that $\textrm{det}(v_1|v_2|\cdots|v_d)=1$.
Let us write $||\cdot||$ the norm on $\mathbb R^d$ given by
$$\left\Vert\sum_{i=1}^dx_iv_i\right\Vert=\max_{i=1,...,d}|x_i| $$
and  $d_0(\cdot,\cdot)$ the metric induced by $||\cdot||$ on $\mathbb R^d$.
Let also $d_1$ be the metric induced by $d_0$ on $\mathbb T^d$.
We define now
$B_u(\delta):=\{y\in E_u\ :\ ||y||\le\delta\}$,
$B_e(\delta):=\{y\in E_e\ :\ ||y||\le\delta\}$
and $B_s(\delta)=\{y\in E_s\ :\ ||y||\le\delta\}$.
Let $|\cdot|$ be the usual euclidean norm on ${\mathbb R}^d$.

Let $r_u$ be the spectral radius of $S^{-1}_{|E_u}$. For every
$\rho_u\in(r_u,1)$, there exists $K>0$ such that, for every integer $n\ge 0$, we have
\begin{equation}\label{hyp1}
\forall h_u\in E_u,\ \
    ||S^n h_u||\ge K\rho_u^{-n}||h_u||
\end{equation}
and
\begin{equation}\label{hyp2}
\forall (h_e,h_s) \in  E_e\times E_s,\ \
    ||S^n (h_e+h_s)||\le K (1+n)^{d_e}||h_e+h_s||.
\end{equation}
Let $\rho_u\in(r_u,1)$ and $K$ satisfying (\ref{hyp1}) and (\ref{hyp2}).
Let $m_u$, $m_e$, $m_s$ be the Lebesgue measure on $E_u$ (in the basis $v_1,...,v_{d_u}$),
$E_e$ (in the basis $v_{d_u+1},...,v_{d_u+d_e}$) and $E_s$ (in the basis
$v_{d_u+d_e+1},...,v_d$) respectively.
Observe that
$d\lambda(h_u+h_e+h_s)=dm_u(h_u)dm_e(h_e)dm_s(h_s)$.

The properties satisfied by the filtration considered in \cite{Lind,SLB}
and enabling the use of a martingale approximation method \`a la Gordin
will be crucial here. Given a finite partition $\mathcal P$ of $\mathbb T^d$, we define
the measurable partition $\mathcal P_0^{\infty}$ by~:
$$\forall \bar x\in\mathbb T^d,\ \
 \mathcal P_0^{\infty}(\bar x):=\bigcap_{k\ge 0}T^k\mathcal P(T^{-k}(\bar x))$$
and, for every integer $n$, the $\sigma$-algebra $\mathcal F_n$ generated by
$$\forall \bar x\in\mathbb T^d,\ \
  \mathcal P_{-n}^{\infty}(\bar x):=\bigcap_{k\ge -n}T^k\mathcal P(T^{-k}(\bar x))
      =T^{-n}(\mathcal P_0^{\infty}(T^n(\bar x)).$$
These definitions coincide with the ones of \cite{SLB} applied to
the ergodic toral automorphism $T^{-1}$.
We obviously have $\mathcal F_n\subseteq\mathcal F_{n+1}=T^{-1}\mathcal F_n$.
Let $r_0>0$ be such that $(h_u,h_e,h_s)\mapsto \overline{h_u+h_e+h_s}$ defines a diffeomorphism
from $B_u(r_0)\times B_e(r_0)\times B_s(r_0)$ on its image in $\mathbb T^d$.
Observe that, for every $\bar x\in\mathbb T^d$,
on the set $\bar x+B_u(r_0)+B_e(r_0)+B_s(r_0)$, we have
$d\bar\lambda(\bar x+\overline{h_u}+\overline{h_e}+\overline{h_s})=dm_u(h_u)dm_e(h_e)dm_s(h_s).$

\begin{Proposition}[\cite{Lind,SLB} applied to $T^{-1}$]\label{partition}
There exist some $Q>0$, $K_0>0$, $\alpha\in(0,1)$
and some finite partition $\mathcal P$ of $\mathbb T^d$
whose elements are of the form
$\sum_{i=1}^dI_i\overline{v_i}$ where the $I_i$ are intervals with diameter smaller than
$\min(r_0,K)$ such that,
for almost every $\bar x\in\mathbb T^d$,
\begin{itemize}
\item[{\bf 1.}] the local leaf $\mathcal P_0^{\infty}(\bar x)$
of $\mathcal P_0^{\infty}$ containing $\bar x$ is a bounded
convex set $\bar x+\overline{F(\bar x)}$, with $0\in F(\bar x)\subseteq E_u$,
$F(\bar x)$ having non-empty interior in $E_u$,
\item[{\bf 2.}] we have
\begin{equation}\label{SLB1}
{\mathbb E}_n(f)(\bar x)=\frac 1{m_u(S^{-n}F(T^n\bar x))}
                 \int_{S^{-n}F(T^n\bar x)}f(\bar x+\overline{h_u})\, dm_u(h_u),
\end{equation}
\item[{\bf 3.}] for every $\gamma>0$, we have
\begin{equation}\label{SLB2}
m_u(\partial (F(\bar x))(\gamma))\le Q\gamma,
\end{equation}
where
$$\partial F(\beta):=\{y\in F\ :\ d(y,\partial F)\le\beta \}, $$
\item[{\bf 4.}] for every $\mathbf k\in\mathbb Z^d\setminus\{0\}$, for every integer $n\ge 0$,
\begin{equation}\label{SLB3}
\left|{\mathbb E}_{-n}(e^{2i\pi\langle \mathbf k,\cdot\rangle})(\bar x) \right|
    \le \frac {K_0}{m_u(F(T^{-n}(x)))} | \mathbf k|^{d_e+d_s}\alpha^n,
\end{equation}
\item[{\bf 5.}] for every $\beta\in(0,1)$,
\begin{equation}\label{SLB4}
\exists L>0,\ \forall n\ge 0,\ \bar\lambda(m_u(F(\cdot))<\beta^n)\le L\beta^{n/ d_u} .
\end{equation}
\end{itemize}
\end{Proposition}
\noindent{\bf Proof.}
The first item comes from Proposition II.1 of \cite{SLB}. Item 2 comes from
the formula given after Lemma II.2 of \cite{SLB}. Item 3 follows from Lemma III.1 of \cite{SLB}
and from the fact that the numbers $a(\mathcal P_0^\infty(\cdot))$ considered in \cite{SLB}
are uniformly bounded. Item 4 comes from Proposition III.3 of \cite{SLB} and from the uniform
boundedness of $a(\mathcal P_0^\infty(\cdot))$.
Item 5 comes from the proof of Proposition II.1 of \cite{SLB}. $\square$
\medskip

\noindent According to the first item of Proposition \ref{partition} and to (\ref{hyp1}),
there exists $c_u>0$ such that, for almost every $\bar x
\in\mathbb T^d$ and every $n\ge 1$, we have
\begin{equation}\label{hyp3}
\sup_{h_u\in S^{-n}F(T^n(\bar x))}\left|h_u\right|\le c_u\rho_u^n.
\end{equation}

\begin{Proposition}\label{terme1}
Let $p\ge 2$ and  $q$ be its conjugate exponent. Let $\theta>0$ and $f:\mathbb T^d\rightarrow \mathbb R$ be a centered function with Fourier coefficients
$(c_{\mathbf k})_{\mathbf k\in\mathbb Z^d}$ satisfying
\beq \label{condfourierq}
\sum_{|\mathbf k|\ge b}|c_{\mathbf k}|^q\le R \log^{-\theta}(b) \, .
\eeq
Then
$$\left\Vert{\mathbb E}_0 (f\circ T^n) \right\Vert_p=
  \left\Vert{\mathbb E}_{-n} (f) \right\Vert_p=O( n^{-\theta(p-1)/p}) \, .$$
\end{Proposition}

\noindent{\bf Proof.}
Recall first that ${\mathbb E}_0 (f\circ T^n)={\mathbb E}_{-n} (f)\circ T^n$.
Let us consider $\alpha$ satisfying (\ref{SLB3}).
Let $\beta:=\alpha^{1/2}$, $\gamma:=\max(\alpha^{ p/2},
  \beta^{1 /d_u} )$ and
$\mathcal V_n:=\left\{\bar x\in\mathbb T^d\ :\ m_u(F(T^{-n}(\bar x))\ge \beta^n\right\}$.
Let $b(n):=\big [ \gamma^{- n/(2p(d+d_e+d_s))} \big ]$. Let us write
\beq \label{decf}
f= f_{1,n} + f_{2,n} \ \text{where $f_{1,n}:=\sum_{|\mathbf k|< b(n)}c_{\mathbf k}e^{2i\pi\langle\mathbf k,\cdot\rangle}
   $ and $
      f_{2,n}:=\sum_{|\mathbf k|\ge b(n)}c_{\mathbf k}e^{2i\pi\langle\mathbf k,\cdot\rangle}$}.
      \eeq
We have
\begin{eqnarray*}
\int_{\mathcal V_n} |{\mathbb E}_{-n} (f_{1,n})|^p\, d\bar\lambda
&\leq & \esup_{\bar x\in\mathcal V_n}
      \Big (\sum_{|\mathbf k|\le b(n)}|c_{\mathbf k}|
          \big |{\mathbb E}_{-n}(e^{2i\pi\langle \mathbf k,\cdot\rangle})(\bar x)
        \big |\Big )^p\\
&\le&  \Big (\sum_{|\mathbf k|\le b(n)}|c_{\mathbf k}|
    {K_0}\beta^{-n}| \mathbf k|^{d_e+d_s}\alpha^n\Big )^p,
\end{eqnarray*}
according to (\ref{SLB3}) and thanks to the definition of $\mathcal V_n$.
Now, since $\beta=\alpha^{ 1/2}$, we get
$$\int_{\mathcal V_n} |{\mathbb E}_{-n} (f_{1,n})|^p\, d\bar\lambda
\le 3^{dp}||f||_1^p K_0^p \alpha^{\frac {np}2} (b(n))^{p(d+d_e+d_s)}.$$
Hence
\begin{equation}\label{p1adapnormp}
\int_{\mathcal V_n} |{\mathbb E}_{-n} (f_{1,n})|^p\, d\bar\lambda
 = O( \gamma^n (b(n))^{p(d+d_e+d_s)}  )=O(\gamma^{n/2}).
\end{equation}
Moreover, thanks to (\ref{SLB4}), we have
\begin{eqnarray}
\int_{\mathcal V_n^c} |{\mathbb E}_{-n} (f_{1,n})|^p\, d\bar\lambda
      &\le& \bar\lambda(\mathcal V_n^c)\Bigl (\sum_{|\mathbf k|\le b(n)}|c_{\mathbf k}|\Bigr
     )^p \nonumber\\
    &=&O((b(n))^{dp}\beta^{n/d_u})=O((b(n))^{dp}\gamma^n)=O(\gamma^{ n/2}).
\label{p2adapnormp}
\end{eqnarray}
Since $p\ge 2$ and since $p/q=p-1$, thanks to (\ref{condfourierq}), we have
\begin{equation} \label{p3adapnormp}
 \Vert{\mathbb E}_{-n} (f_{2,n})\Vert_p^p
      \le \Vert f_{2,n}\Vert_p^p
       \le \Big(\sum_{|\mathbf k|\ge b(n)}|c_{\mathbf k}|^q\Big)^{ p/q}
       \le R^{p-1}(\log(b(n)))^{-\theta(p-1)} \ll  n^{-\theta(p-1)} \, .
\end{equation}
Combining \eqref{p1adapnormp}, \eqref{p2adapnormp} and \eqref{p3adapnormp},
the proposition follows. $\square$

\begin{Proposition}\label{terme3}
Under the assumptions of Proposition \ref{terme1},
$$\left\Vert{\mathbb E}_0 (f\circ T^{-n})-f\right\Vert_{p}=
  \left\Vert{\mathbb E}_n (f)-f\right\Vert_{p}=O( n^{-\theta(p-1)/p})\, . $$
\end{Proposition}
\noindent{\bf Proof.} We consider the decomposition (\ref{decf}) with $b(n)$ defined by
$b(n)= \big  [ \rho_u^{- n/(2(d+1))}  \big ] $. We have
\begin{eqnarray*}
\Vert{\mathbb E}_n (f_{1,n})-f_{1,n}\Vert_{p}
  &\le&   \Vert{\mathbb E} _n(f_{1,n})-f_{1,n} \Vert_{\infty}\\
  &\le& \sum_{|\mathbf k|\le b(n)}|c_{\mathbf k}|
      \Vert \E_n ( e^{2i\pi\langle\mathbf k,\cdot \rangle})- e^{2i\pi\langle\mathbf k,\cdot \rangle}\Vert_{\infty}\\
  &\le&  \sum_{|\mathbf k|\le b(n)}|c_{\mathbf k}|
      2\pi |\mathbf k| c_u\rho_u^{n} \, ,
\end{eqnarray*}
according to (\ref{SLB1}) and to (\ref{hyp3}).
Therefore
\begin{equation} \label{p1normnadapt}
\left\Vert{\mathbb E}_n (f_{1,n})-f_{1,n}\right\Vert_{p}
     \ll (b(n))^{d+1}\rho_u^{n}\ll \rho_u^{ n/2} .
\end{equation}
Moreover, thanks to (\ref{condfourierq}), we have
\begin{align} \label{p2normnadapt}
\Vert{\mathbb E}_n (f_{2,n})-f_{2,n}\Vert_{p}^p
   &  \le 2^p \Vert f_{2,n}\Vert_p^p\le 2^p
       \Big (\sum_{|\mathbf k|\ge b(n)}|c_{\mathbf k}|^q\Big )^{ p/q} \nonumber \\
      &  \le 2^pR^{p-1}(\log(b(n)))^{-\theta(p-1)}\ll  n^{-\theta(p-1)} \, .
\end{align}
Considering \eqref{p1normnadapt} and \eqref{p2normnadapt}, the proposition follows. $\square$

\begin{Proposition}\label{terme2}
Let $p \in [2,4]$ and set
$S_n(f) := \sum_{k=1}^n f\circ T^k$ with $f:\mathbb T^d\rightarrow \mathbb R$ be a
centered function with Fourier coefficients satisfying \eqref{condfourierq} with $\theta>0$
and
\beq \label{condfourier2}
\sum_{|\mathbf k|\ge b}|c_{\mathbf k}|^2\le R \log^{-\beta }(b)  \ \text{ for some $\beta>1$}\, .
\eeq
Set
\beq \label{definitiondec}
m:= \Big [ - \frac{4(d_e+d_s)\log(r)} {\log (\alpha)}  \Big ] +1 \, .
\eeq
where   $r$ is the spectral radius of $S$. Then
$$
\Vert \E_{-nm }(S^2_n(f)) - \E(S_n^2(f)) \Vert_{p/2} \ll n^{2 -2 \theta (p-1)/p} +  n^{(3- \beta)/2} \, .
$$
\end{Proposition}
\noindent{\bf Proof.} Let $\beta:=\alpha^{ 1/2}$,
$\mathcal V_{nm}:=\left\{\bar x\in\mathbb T^d\ :\ m_u(F(T^{-nm}(\bar x))\ge \beta^{nm}
  \right\}$,
$\gamma:=\max(\alpha^{ p/8},
  \beta^{ 1/d_u})$ and
\beq \label{choixbndouble} b(n):=\Big [ \gamma^{ {n \, m }/(p(2d+d_e+d_s))}  \Big ] \, . \eeq
We consider the decomposition (\ref{decf}) with $b(n)$ defined by (\ref{choixbndouble}) and we set
$$
S_{1,n}(f):=\sum_{k=1}^n f_{1,n}\circ T^k \ \text{ and } S_{2,n}(f):=\sum_{k=1}^n f_{2,n}\circ T^k \, .
$$
First, we note that
\begin{multline*}
 \|\E_{-nm}(S^2_n(f))  - \E(S_n^2(f)) \Vert_{p/2} \leq \Vert \E_{-nm}(S^2_{1,n}(f)) -
 \E(S^2_{1,n}(f)) \Vert_{p/2} \\
 + \Vert \E_{-nm}(S^2_{2,n}(f)) - \E(S^2_{2,n}(f)) \Vert_{p/2} +2 \Vert \E_{-nm}(S_{1,n}(f)S_{2,n}(f)) - \E(S_{1,n}(f)S_{2,n}(f)) \Vert_{p/2} \\
\leq \Vert \E_{-nm}(S^2_{1,n}(f)) - \E(S^2_{1,n}(f)) \Vert_{p/2} + 2 \Vert S_{2,n}(f) \Vert_{p}^2 + 4 \Vert {\mathbb E}_{-nm}(S_{1,n}(f)S_{2,n}(f))  \Vert_{p/2} \, .
\end{multline*}
Next using \eqref{fact}, we get that
\begin{multline*}
\Vert \E_{-nm}(S_{1,n}(f)  S_{2,n}(f))  \Vert_{p/2} \leq \Vert \E_{-nm}(S^2_{1,n}(f)) - \E(S^2_{1,n}(f)) \Vert_{p/2} \\
  + \Vert S_{2,n}(f)  \Vert^2_{p} + \Vert S_{1,n}(f)  \Vert_{2} \Vert S_{2,n}(f) \Vert_{2} \\
 \leq \Vert \E_{-nm}(S^2_{1,n}(f)) - \E(S^2_{1,n}(f)) \Vert_{p/2} +
 2 \Vert S_{2,n}(f)  \Vert^2_{p} + \Vert S_{n}(f)  \Vert_{2} \Vert S_{2,n}(f) \Vert_{2} \, .
\end{multline*}
By Propositions \ref{terme1} and \ref{terme3}, \eqref{condfourier2} implies that
$$
\sum_{n >0} \frac{\Vert {\mathbb E}_{-n} (f) \Vert_2}{n^{1/2}} < \infty \ \text{ and } \
  \sum_{n >0} \frac{\Vert f - {\mathbb E}_n (f) \Vert_2}{n^{1/2}} < \infty \, ,
$$
which  yields \eqref{Gornorm} with $p=2$, and then $\Vert S_{n}(f)  \Vert_{2} \ll \sqrt{n}$. Therefore, we get overall that
\begin{multline} \label{doubp1}
\Vert  \E_{-nm}(S^2_n(f))  - \E(S_n^2(f)) \Vert_{p/2} \ll \Vert \E_{-nm}(S^2_{1,n}(f)) - \E(S^2_{1,n}(f)) \Vert_{p/2} \\
+ \Vert S_{2,n}(f)  \Vert^2_{p} +  \sqrt{n} \Vert S_{2,n}(f)  \Vert_{2}\, .
\end{multline}
Since $p\ge 2$ and $p/q=p-1$, \eqref{condfourierq} implies that
\begin{align}\label{doubp2}
\Vert S_{2,n}(f)  \Vert_{p} & \leq n \Vert f_{2,n} \Vert_{p}
     \le n \Big ( \sum_{|\mathbf k|\ge b(n)}|c_{\mathbf k}|^q \Big )^{ 1/q} \nonumber \\
           & \le n R^{(p-1)/p}(\log(b(n)))^{-\theta(p-1)/p}\ll  n^{1-\theta(p-1)/p} \, .
              \end{align}
Similarly using \eqref{condfourier2}, we get that
\begin{equation}\label{doubp3}
\Vert S_{2,n}(f)  \Vert_{2} \leq n \Vert f_{2,n} \Vert_{2}
     \ll n^{1 -\beta/2} \, .
              \end{equation}
We deal now with the first term in the right hand side of \eqref{doubp1}. With this aim,
we first observe that, for any non negative integer $\ell$,
$ e^{2i\pi\langle\mathbf k,T^\ell(\cdot)\rangle}=
       e^{2i\pi\langle{}^tS^\ell\mathbf k,\cdot\rangle},$
where ${}^tS^\ell$ is the transposed matrix of $S^{\ell}$. Therefore,
\begin{align*}
\int_{\mathcal V_{ n m }} \big |{\mathbb E}_{-nm} (f_{1,n}.f_{1,n}& \circ T^\ell)
         - {\mathbb E} (f_{1,n}.f_{1,n}\circ T^\ell )   \big |^{p/2}\,
         d\bar\lambda \\
&\le \esup_{\bar x\in\mathcal V_{ n m }}\Big (\sum_{|\mathbf k|,|\mathbf m|\le b(n):
                     \mathbf k+{}^tS^\ell\mathbf m\ne 0}|c_{\mathbf k}|
           |c_{\mathbf m}|
          \big |{\mathbb E}_{-nm}(e^{2i\pi\langle \mathbf k+{}^tS^\ell\mathbf m
      ,\cdot\rangle})(\bar x) \big |\Big )^{p/2}\\
&\le  \Big (\sum_{|\mathbf k|,|\mathbf m|\le b(n)}
      |c_{\mathbf k}||c_{\mathbf m}|
       {K_0}\beta^{-n m}| \mathbf k+{}^tS^\ell\mathbf m|^{d_e+d_s}\alpha^{n m}\Big )^{p/2},
\end{align*}
according to (\ref{SLB3}) and to the definition of $\mathcal V_{nm}$.
It follows that
\begin{align*}
\int_{\mathcal V_{ n m }} \big |{\mathbb E}_{-nm} (f_{1,n}.f_{1,n}& \circ T^\ell)
         - {\mathbb E} (f_{1,n}.f_{1,n}\circ T^\ell )   \big |^{p/2}\,
         d\bar\lambda  \\
&\le \Big ( \sum_{|\mathbf k|,|\mathbf m|\le b(n)}
      \Vert f \Vert_1^2 {K_0}\beta^{-n m } (| \mathbf k|+ r^\ell|\mathbf m|)^{d_e+d_s}\alpha^{n m}
     \Big )^{p/2}\\
      &\ll \alpha^{\frac {n m p}4} r^{p\ell(d_e+d_s)/2}
     (b(n))^{p(2d+d_e+d_s)/2} \, .
\end{align*}
Hence, since $\gamma \geq \alpha^{p/8}$, $m \geq  {4(d_e+d_s)\log(r)} /{\log (1/\alpha)}$,
and according to the definition of $b(n)$, we have
\begin{equation} \label{p1normp2} \sup_{\ell\in \{0,\ldots, n \}}\int_{\mathcal V_{ n m}}
\Big|{\mathbb E} _{-nm}(f_{1,n}.f_{1,n}\circ T^\ell)-
   {\mathbb E} (f_{1,n}.f_{1,n}\circ T^\ell)\Big|
     ^{p/2}\, d\bar\lambda
  \ll  \alpha^{ { 3 n m p}/{16}} r^{p n (d_e+d_s)/2}  \\
  \ll \gamma^{{  n m }/{2}} .
\end{equation}
Moreover, for any non negative integer $ \ell$,
\begin{multline} \label{p2normp2}
\int_{\mathcal V_{n m}^c} \Big|{\mathbb E}_{-nm}  (f_{1,n}.f_{1,n}\circ T^\ell)
    \Big|^{p/2}\,   d\bar\lambda
      \le \bar\lambda(\mathcal V_{n m}^c)\Bigl (\sum_{|\mathbf k|,|\mathbf m|\le b(n)}
     |c_{\mathbf k}||c_{\mathbf m}|\Bigr )^{p/2}  \\
 \ll(b(n))^{dp}\beta^{n m /d_u} \ll (b(n))^{dp}\gamma^{n m} \ll \gamma^{ {  n m }/{2}}  \, ,
\end{multline}
according to (\ref{SLB4}) and to the definition of $b(n)$ and of $\gamma$.
Combining \eqref{p1normp2}  and \eqref{p2normp2}, we then derive that
\begin{align} \label{p3normp2}
\Vert \E_{-nm}& (S^2_{1,n}(f))  - \E(S^2_{1,n}(f)) \Vert_{p/2} \nonumber \\
& \leq 2 \sum_{i=1}^n \sum_{j=0}^{n-i}\Vert \E_{-nm}(f_{1,n}\circ T^i f_{1,n}\circ T^{i+j} )- \E(f_{1,n}\circ T^i f_{1,n}\circ T^{i+j} ) \Vert_{p/2} \nonumber \\
& \leq  n^2\sup_{\ell \in \{0,\ldots, n \}}\Vert \E_{-nm}(f_{1,n} f_{1,n}\circ T^{\ell} )- \E(f_{1,n} f_{1,n}\circ T^{\ell} ) \Vert_{p/2} \ll n^2\gamma^{ {  n m }/{p}}  \, .
\end{align}
Considering \eqref{doubp2}, \eqref{doubp3} and  \eqref{p3normp2} in \eqref{doubp1}, the proposition follows. $\square$

\subsection{End of the proof of  Theorem  \ref{auto-tore2}}
Propositions \ref{terme1} and \ref{terme3} give (\ref{Cond1cob*}) provided
(\ref{condF1}) is satisfied.
Propositions  \ref{terme1} and \ref{terme3} give (\ref{Cond1cob**})
(when $p\in]2,4]$) and (\ref{Cond1cob**p=4}) (when $p=4$),
provided (\ref{condF2}) is satisfied.
Finally, Proposition \ref{terme2} gives (\ref{condcarre}) provided
(\ref{condF1}) and (\ref{condF2}) are satisfied.
The proof follows now from
Theorem \ref{ThNAS} when $p\in ]2,4[$ and from Theorem \ref{ThNASp=4} when $p=4$. $\square$

\section{Appendix}
As in Section \ref{MR}, let $P_k(X)= \E_{k}(X)-\E_{k-1}(X)$.
\setcounter{equation}{0}
\begin{Lemma} \label{complemmap0xi} Let $p\in [2, \infty[$.
Then, for any real $1 \leq q \leq p$ and any positive integer $n$,
$$
\sum_{k \geq 2 n}\Vert P_0(X_{k}) \Vert^q_p \ll \sum_{k \geq n  } \frac{\Vert \E_0(X_k)\Vert^q_p}{k^{q/p}} \, \text{ and } \sum_{k \geq 2 n}\Vert P_0(X_{-k}) \Vert^q_p \ll \sum_{k \geq n  } \frac{\Vert X_{-k} - \E_0(X_{-k} )\Vert^q_p}{k^{q/p}} \, .
$$
\end{Lemma}
\noindent {\bf Proof.} The first inequality is Lemma 5.1 in \cite{DDM}.
To prove the second one, we first consider the case  $p>q$ and we
follow the lines of the proof Lemma 5.1 in \cite{DDM} with $P_k(X_{0})$ replacing
$P_{-k} (X_0)$.
We get that
$$\sum_{k\ge 2n}\Vert P_0(X_{-k})\Vert^q_p \ll \sum_{k\ge n+1}k^{-\frac qp}
      \Big(\sum_{\ell\ge k}\Vert P_0(X_{-\ell})\Vert^p_p\Big)^{ q/p}.$$
Now, we notice that, by the Rosenthal's inequality given in Theorem 2.12 of \cite{HallHeyde},
there
exists a constant $c_{p}$ depending only on $p$ such that
\begin{multline} \label{ros1}
\sum_{\ell\ge k}\Vert P_{0}(X_{-\ell}) \Vert _{p}^{p}
  = \sum_{\ell  \geq k } \Vert P_{\ell}(X_{0}) \Vert _{p}^{p} \\
 \leq  c_p
   \Big \Vert \sum_{\ell  \geq k } P_{\ell}(X_{0}) \Big \Vert _{p}^{p}=c_p \Vert X_0 - \E_k (X_0) \Vert _{p}^{p} =c_p \Vert X_{-k} - \E_0 (X_{-k}) \Vert _{p}^{p} \, .
\end{multline}
Now when $p=q$, inequality \eqref{ros1} together with the fact that by Claim
\ref{claimfacile}, for any integer $k$ in $[n+1,2n]$,
$ \Vert X_0 - \E_{2n} (X_0) \Vert _{p}^{p} \leq 2^p \Vert X_0 - \E_{k} (X_0) \Vert _{p}^{p}$
imply the result.
Indeed we have
$$\sum_{k\ge 2n} \Vert P_{0}(X_{-\ell}) \Vert _{p}^{p}\le c_p
     \Vert X_0 - \E_{2n} (X_0) \Vert _{p}^{p}\ll \sum_{k=n+1}^{2n}k^{-1}
  \Vert X_0 - \E_{k} (X_0) \Vert _{p}^{p} \, .\quad \square$$

\bigskip
\begin{Proposition} \label{approxrnq} Let  $p\in [1, \infty[$ and assume that
\begin{equation}\label{Gor}
  \text{the series} \quad d_0=\sum_{i \in {\mathbb Z}} P_0(X_i) \quad
  \text{converges in ${\mathbb L}^p$.}
\end{equation}
Let $M_n:= \sum_{i=1}^n d_0 \circ T^i $ and $R_n:=S_n-M_n$. Then, for any positive integers $n$ and  $N$,
\begin{align*}
 R_n  & =   \E_0 (S_n ) -\E_0 (S_{N} ) \circ T^n + \E_{-n} (S_{N}) \circ T^n -  \sum_{k=1}^n \sum_{j \geq n+N +1} P_k(X_j) \\
& + S_n - {\mathbb E}_n (S_n ) - (  \E_{n+N} (  S_{N} - \E_N (S_{N} ) ) \circ T^{-N}  - \sum_{k=1}^n \sum_{j\geq N} P_k(X_{-j}) \, ,
\end{align*}
and
\begin{align*}
\Vert R_n \Vert^{p'}_{p} \ll   \Vert \E_0 (S_{n}  ) \Vert_p^{p'} & + \Vert \E_0 (S_{N} ) \Vert_p^{p'} +    \Vert S_n - {\mathbb E}_n (S_n )  \Vert^{p'}_p + \Vert S_{N} - \E_N (S_{N} ) \Vert^{p'}_p \\
& + \sum_{k=1}^n
\big \Vert  \sum_{j \geq k+N}  P_0(X_j) \big \Vert_p^{p'} + \sum_{k=1}^n
\big \Vert  \sum_{j \geq k+N}  P_0(X_{-j}) \big \Vert_p^{p'} \, ,
\end{align*}
where $p'=\min ( 2,p)$.
\end{Proposition}

\noindent{\bf Proof of Proposition \ref{approxrnq}.} Notice first that the following decomposition is valid:  for any positive integer $n$,
\begin{equation} \label{decRn*}
R_n   =  \sum_{k=1}^n \Big ( X_k - \sum_{j =1}^n P_j(X_k) \Big )
  - \sum_{k=1}^n \sum_{j \geq n+1} P_k(X_j) - \sum_{k=1}^n \sum_{j=0}^{\infty} P_k(X_{-j})
  =  R_{n,1} + R_{n,2} \, ,
\end{equation}
where
\beq \label{defrn12}
R_{n,1} :=  {\mathbb E}_0 (S_n) - \sum_{k=1}^n \sum_{j \geq n+1} P_k(X_j)
\text{ and }R_{n,2} := S_n - {\mathbb E}_n (S_n )  -
\sum_{k=1}^n \sum_{j=0}^{\infty} P_k(X_{-j}) \, .
\eeq
Let $N$ be a positive integer. According to item 1 of Proposition 2.1 in \cite{DDM},
\beq \label{brn1*}
 R_{n,1} =  \E_0 (S_n) -\E_n (S_{n+N} - S_n ) +  \E_0 (S_{n+N} - S_n ) -
    \sum_{k=1}^n \sum_{j \geq n+N +1} P_k(X_j)\, .
\eeq
On an other hand, we write that  $
\sum_{j=0}^{\infty} P_k(X_{-j}) = \sum_{j=0}^{N -1} P_k(X_{-j}) + \sum_{j\geq N} P_k(X_{-j})
$. Therefore
\beq \label{normeprn2*}
 R_{n,2} =S_n - {\mathbb E}_n (S_n) - (  \E_{n+N} (  S_{N} - \E_N (S_{N} ) ) \circ T^{-N}  - \sum_{k=1}^n \sum_{j\geq N} P_k(X_{-j})  \, .
\eeq
Starting from \eqref{decRn*} and considering \eqref{brn1*} and \eqref{normeprn2*},
the first part follows. We turn now to the second part of the proposition. Applying
Burkholder's inequality and using  stationarity, we obtain that
there exists a positive constant $c_p$ such that, for any positive integer $n$,
\begin{eqnarray} \label{Burkh1}
\Big \Vert  \sum_{k=1}^n \sum_{j \geq n+N +1} P_k(X_j) \Big \Vert^{p'}_{p}   \leq   c_p
  \sum_{k=1}^n \Big \Vert \sum_{j \geq n+N +1} P_k(X_j) \Big \Vert_p^{p'}  = c_p
   \sum_{k=1}^n  \Big \Vert  \sum_{j \geq N+k}  P_0(X_j) \Big \Vert_p^{p'} \, ,
\end{eqnarray}
and
\begin{eqnarray} \label{Burkh2}
\Big \Vert  \sum_{k=1}^n \sum_{j\geq N} P_k(X_{-j})  \Big \Vert^{p'}_{p}   \leq
    c_p  \sum_{k=1}^n \Big \Vert \sum_{j\geq N} P_k(X_{-j})  \Big \Vert_p^{p'}  = c_p
    \sum_{k=1}^n  \Big \Vert  \sum_{j \geq N+k}  P_0(X_{-j}) \Big \Vert_p^{p'} \, .
\end{eqnarray}
The second part of the proposition  follows from item 1 by taking into account  stationarity and by
considering the bounds \eqref{Burkh1} and \eqref{Burkh2}. $\square$

\end{document}